# Multi-agent Modeling and Optimal Pumping Control of Magnetic Artificial Cilia
# (CCSICC)

Shuangshuang Yu[1,2], Zheng Ning[2], and Ge Chen[2]

[1)] School of Mathematical Sciences, University of Chinese Academy of Sciences, Beijing 100049, China
(yushaungshuang@amss.ac.cn)

[2)] The Key Laboratory of Systems and Control, Academy of Mathematics and Systems Science, Chinese Academy of Sciences, Beijing 100190, China (ningzheng@amss.ac.cn, chenge@amss.ac.cn)

*Abstract*—Tiny cilia drive the flow of surrounding fluids through asymmetric jumping, which is one of the main ways for biological organisms to control fluid transport at the microscale. Due to its huge application prospects in medical and environmental treatment fields, artificial cilia have attracted widespread research interest in recent years. However, how to model and optimize artificial cilia is currently a common challenge faced by scholars. We model a single artificial cilium driven by a magnetic field as a multi-agent system, where each agent is a magnetic bead, and the interactions between beads are influenced by the magnetic field. Our system is driven by controlling the magnetic field input to achieve fluid transport at low Reynolds number. In order to quantify the flow conveying capacity, we introduce the pumping performance and propose an optimal control problem for pumping performance，and then give its numerical solution. The calculation results indicate that our model and optimal control algorithm can significantly improve the pumping performance of a single cilia.

# 磁性人工纤毛的多个体建模与泵送最优控制

郁霜霜 [1,2]    宁政 [2]    陈鸽 [2*]

[1)]中国科学院大学数学科学学院，北京 100049，中国
[2)] 中国科学院数学与系统科学研究院系统控制重点实验室，北京 100190，中国

摘 要 细小的纤毛通过不对称跳动来带动周围液体流动，是生物有机体在微尺度上控制流体输送的主要方式之一。由于在医疗和环境处理等领域的巨大应用前景，人工纤毛近年来吸引了广泛研究兴趣。然而，如何对人工纤毛建模与优化是目前学者们所面临的共同难题。我们将受磁场驱动的单根人工纤毛建模成一个多个体系统，其中每个个体为一个磁性微粒，个体间相互作用受磁场影响。我们通过控制磁场输入来驱动系统，实现低雷诺数下的流体输送。为了衡量流量输送能力，我们引入泵送性能这一评价指标，并提出了泵送性能的最优控制问题及其数值解法。计算结果表明，我们的模型和最优控制算法可能显著改善单根纤毛的泵送性能。

关键词 微纳米机器人，人工纤毛，泵送性能，多个体系统，最优控制

## 1．引言

近年来，已经有许多关于血液循环系统的仿真研究[1-2]。血液作为循环流动在血管与心脏中的一种运输组织，它将氧气运送到各个器官并带离细胞代谢产生的废物，这意味着如果血液循环出现问题，那么可能产生血栓，引起

心肌梗塞、脑中风、脑血管溢血等严重疾病。正常情况下生物体中血液雷诺数比较低，因此本文将从微观角度研究低雷诺数下的流体输送最优控制问题。

在低雷诺数下的微流体中，惯性力不显著，而粘性效应占据主导，这意味着往复运动形成的流体置换并不会产生净输送量。而纤毛作为许多生物体的运动器官，其独特的运动模式能够适用于低雷诺数下的流体输送。人工纤毛作为近年来最受关注的微纳机器人之一，其驱动方式可分为外部物理场驱动（通常有磁场、光场、声场、电场等）、化学燃料驱动以及这些驱动方式的混合驱动[3-5]。其中，磁性驱动由于其对生物无害、远程可控、快速响应、精确控制、低能耗等显著优点，在靶向治疗、微创外科手术、体内可视化检测、重金属检测、污染物降解等领域有非常重要的潜在应用[6-9]，因此引起了国内外研究人员的广泛兴趣。

目前已经存在很多关于纤毛泵送微流体的研究工作，然而当下的研究大多是基于实验层面上的，例如 Hanasoge 等[10]利用光刻技术制造人工纤毛，并使用磁场驱动实现在≈1 Pa 的压强下每分钟 11 毫升的泵送流速；Toonder 等[11]制造了由聚合物结构组成的人造纤毛，并将其集成到微流体通道中进行流体泵送与操纵；Vilfan 等[12]创建了自组装人工纤毛的有序阵列，并证明了它们可以用于泵送微流体。另一方面，人工纤毛的理论建模对于理解实验结果和指导未来的实验非常重要。现已有一些研究给出了二维/三维人工纤毛的运动建模，例如 Zhang 等[13]对人工纤毛的异时搏动行为进行了数值分析，并提出了一个双纤毛模型来实现更高的泵送效率；Gilpin 等[14]总结了单纤毛的内部工作原理和流体动力学模型，以及多纤毛系统在主动滤波、噪声鲁棒性和异时波等方面呈现出的涌现行为；Khaderi 等[15]建立了一种三维 FSI 模型，用于研究磁性人工纤毛与流体间的相互作用，并利用该模型对磁性人工纤毛产生的流体流动进行了模拟研究与分析。但目前为止，人工纤毛泵送微流体能力的最优控制理论研究基本空白。

由于低雷诺数下纤毛泵产生的净输送量依赖于纤毛的不对称跳动，而这种跳动可以通过外部磁场的频率和相位输入进行控制[16-17]。一个自然的问题是：怎样的磁场控制能够诱导纤毛产生最强泵送量？本文拟对微流体中低雷诺数下的纤毛泵送性能进行数学建模与分析，由于人工纤毛由具有相互作用的多个磁性微粒组成，因此首先将其建模成一个受磁场驱动的多个体系统，然后研究该系统泵送性能的最优磁场控制理论与算法。

本文组织结构如下：第 2 节首先给出了对人工纤毛的多个体建模，然后引入了微流体输送能力的衡量指标—泵送性能；第 3 节介绍了人工纤毛泵送性能的最优控制问题与相关理论；第 4 节利用共轭梯度法对最优控制问题给出了数值解；第 5 节总结了该论文并给出未来若干研究问题。

## 2．人工纤毛的多个体建模与泵送性能

### 2.1 人工纤毛的多个体建模

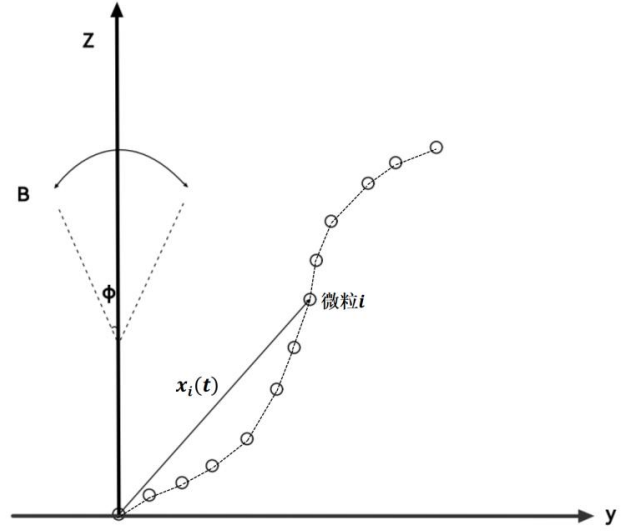

图 1　某个$t$时刻人工纤毛长链状态示意图

如图 1 所示，我们考虑的人工纤毛是由$n$个具有相互作用的磁性微粒构成的长链，其在微流体中进行平面运动且一端固定。个体$i$在$t$时刻的坐标向量记作$x_i = x_i(t) = \begin{pmatrix} y_i(t) \\ z_i(t) \end{pmatrix}$。根据文献[17]，长链的总自由能由磁偶极相互作用能、拉伸自由能和弯曲自由能三部分组成，即

$$E = E(t) = E^S + E^B + E^D, \tag{1}$$

其中$E^S = E^S(t)$是$t$时刻的总拉伸自由能，$E^B = E^B(t)$是$t$时刻的弯曲自由能，$E^D = E^D(t)$是$t$时刻的微粒之间的总磁偶极相互作用能。它们的表达式分别如下：

$$E^S = \frac{1}{2} k \sum_{i=2}^{n}(l_i - l_0)^2, \tag{1.1}$$

$$E^B = \frac{A}{l_0} \sum_{i=2}^{n}\bigl(1 - \hat{l}_i \cdot \hat{l}_{i+1}\bigr)^2, \tag{1.2}$$

$$E^D = \frac{4\pi a^6}{9\mu_0} (\chi B)^2 \sum_{i<j} \frac{1-3(\hat{p}\cdot\hat{x}_{ij})^2}{x_{ij}^3}, \tag{1.3}$$

其中$k > 0$为拉伸系数，$A > 0$为抗弯刚度，$a > 0$为磁性微粒半径，真空磁导率$\mu_0 = 4\pi \times 10^{-7} N/A^2$，$\chi \in \mathbb{R}$为磁性微粒的磁化率，$B \geq 0$为磁场强度，$\hat{p}$为磁场方向单位向量；$l_i = l_i(t) = |x_i(t) - x_{i-1}(t)|$，$\hat{l}_i = \hat{l}_i(t) = \frac{x_i(t)-x_{i-1}(t)}{l_i}$，$l_0 = l_0(t) = \frac{1}{n}\sum_{i=1}^{n} l_i(t)$，$x_{ij} = x_{ij}(t) = |x_j(t) - x_i(t)|$，$\hat{x}_{ij} = \hat{x}_{ij}(t) = \frac{x_j(t)-x_i(t)}{x_{ij}}$。

在低雷诺数下，粒子浸入粘性流体(如水)中的运动速度遵循下面的方程[18]：

$$v_i = v_i(t) = \sum_j \boldsymbol{\mu}_{ij} F_j, \quad (2)$$

其中$\boldsymbol{\mu}_{ij} = \boldsymbol{\mu}_{ij}(t)$为$t$时刻的磁性微粒$i$与磁性微粒$j$间的流动性矩阵，$F_j = F_j(t)$为$t$时刻磁性微粒$j$所受到的力。我们采用 Rotne-Prager 方法来近似$\boldsymbol{\mu}_{ij}$ [17]，其自流动性和交叉流动性可分别表示为

$$\boldsymbol{\mu}_{ii} = \frac{1}{6\pi\eta a}\boldsymbol{I}, \quad (3)$$

$$\boldsymbol{\mu}_{ij} = \frac{1}{6\pi\eta x_{ij}}\left[\frac{3}{4}\left(\boldsymbol{I} + \hat{x}_{ij}\hat{x}_{ij}^T\right) + \frac{a^2}{2x_{ij}^2}\left(\boldsymbol{I} - 3\hat{x}_{ij}\hat{x}_{ij}^T\right)\right], \quad (4)$$

其中$\eta > 0$为流体粘度。$F_j$可由下式推出：

$$F_j = -\nabla_{x_j}(E^S + E^B + E^D)。 \quad (5)$$

计算可得

$$F_j^S := -\nabla_{x_j}E^S = -k(l_j - l_0)\hat{l}_j + k(l_{j+1} - l_0)\hat{l}_{j+1}, \quad (6)$$

$$F_j^B := -\nabla_{x_j}E^B$$
$$= \frac{A}{l_0}\left[\frac{l_{j-1}}{l_j} - \left(\frac{l_j \cdot l_{j-1}}{l_j} + \frac{1}{l_{j+1}} + \frac{l_j \cdot l_{j+1}}{l_j}\right)\hat{l}_j\right.$$
$$\left.+ \left(\frac{l_j \cdot l_{j+1}}{l_{j+1}} + \frac{1}{l_j} + \frac{l_{j+1} \cdot l_{j+2}}{l_{j+1}}\right)\hat{l}_{j+1} - \frac{l_{j+2}}{l_{j+1}}\right], \quad (7)$$

$$F_j^D := -\nabla_{x_j}E^D$$
$$= \frac{4\pi a^6}{3\mu_0}(\chi B)^2 \sum_{i>j} \frac{2(\hat{p} \cdot \hat{x}_{ij})\hat{p} + [1 - 5(\hat{p} \cdot \hat{x}_{ij})^2]\hat{x}_{ij}}{x_{ij}^4}。 \quad (8)$$

### 2.2 人工纤毛的泵送性能

本文拟研究在平面上跳动的磁链所引发的流体输送性能。根据文献[18]，积分流由横向平均布莱克张量所决定：

$$\bar{G}(z,z') = \iint G_{yy}^{Blake}(x,y,z,z')\,dx\,dy = \frac{\min(z,z')}{\eta}。 \quad (9)$$

由于流体输送是沿着图 1 中$y$轴方向的，根据[17]整条磁链产生的积分流可以近似为

$$\mathcal{F}(z) = \sum_{i=1}^n \bar{G}(z,z_i)F_{y_i}, \quad (10)$$

其中$z_i$为第$i$个磁珠与$y$轴的垂直距离。

当$z$大于磁链长度$L$时，$z > z_i$对于任意的$i = 1, \cdots, n$都成立，那么$\bar{G}(z,z_i) = \frac{\min(z,z_i)}{\eta} = \frac{z_i}{\eta}$，进而对任意给定的$t$可得$\mathcal{F}(z > L) = \sum_{i=1}^n \frac{z_i(t)}{\eta}F_{y_i}(t)$，它不依赖于$z$。因此考虑在一个跳动周期内$\mathcal{F}(z > L)$的平均值作为流体输送的度量，称之为单根人工纤毛的泵送性能：

$$\int_t^{t+T} \frac{\mathcal{F}(z>L)}{T} dt = \frac{1}{T\eta}\int_t^{t+T} \sum_{i=1}^n z_i(t)F_{y_i}(t)dt。 \quad (11)$$

### 3. 泵送性能的最优控制

为了探究怎样的磁场控制可以诱导实现泵送性能的最优化，我们固定磁场强度，视磁场偏离$z$轴的角度$\phi(t)$为控制变量，即$u = \phi(t)$，$t \in (t_0, t_f)$，$u \in U = [0, 2\pi]$；用磁性微粒的位置函数做状态变量，即$x_i(t) = \begin{pmatrix} y_i(t) \\ z_i(t) \end{pmatrix}$，$i = 1, \cdots, n$，$x(t) = \begin{pmatrix} x_1(t) \\ x_2(t) \\ \vdots \\ x_n(t) \end{pmatrix}$，$x(t)$为$2n$维向量值函数。

不妨认为人工纤毛的初始位置位于$z$轴上，结合(2)式，状态方程及初始条件为：

$$\begin{cases} \dot{x}_i = \sum_j \boldsymbol{\mu}_{ij}F_j, & i = 1, \cdots, n, \\ x(t_0) = (0, l_0, 0, 2l_0, \cdots, 0, nl_0)^T \triangleq x_0。 \end{cases} \quad (12)$$

由第 2 节讨论内容我们知道，本文最优化的性能指标是 Lagrange 型积分性能。由于$F_{y_i} = -\frac{\partial E}{\partial y_i}$，因此最大化单根人工纤毛的泵送性能相当于：

$$\min_{u \in U} J(u) = \min_{u \in U} \frac{1}{T\eta}\int_t^{t+T}\sum_{i=1}^n z_i \frac{\partial E}{\partial y_i} dt。 \quad (13)$$

下面我们将使用文献[19]中提出的极大值原理来解决上述无约束最优控制问题。系统(12)和性能指标$J(u)$对应的 Halmiton 函数为：

$$H(x(t), u(t), \lambda(t)) = -l(x(t), u(t)) + \lambda^T(t)f(x(t), u(t)),$$

其中$l(x(t), u(t)) = \frac{1}{T\eta}\sum_{i=1}^n z_i(t)\nabla_{y_i}\left(E(x(t), u(t))\right)$，

$f(x(t), u(t)) = \begin{pmatrix} \dot{x}_1(t) \\ \dot{x}_2(t) \\ \vdots \\ \dot{x}_n(t) \end{pmatrix}$，$\lambda(t)$是 Lagrange 乘子向量值函数，它为下终值问题：

$$\begin{cases} \dot{\lambda}(t) = -\frac{\partial H}{\partial x}, \\ \lambda(t_f) = 0。 \end{cases} \quad (14)$$

为了表达方便，我们记$\boldsymbol{\mu}_{ij}(t) = \begin{pmatrix} a_{ij}(t) & b_{ij}(t) \\ c_{ij}(t) & d_{ij}(t) \end{pmatrix}$，$i, j = 1, \cdots, 2n$，结合式(3)–(11)，式(12)可写为

$$\begin{cases} \dot{\lambda}(t) = \boldsymbol{S}(t)\lambda(t) + g(t), \\ \lambda(t_f) = 0, \end{cases} \quad (15)$$

其中$\boldsymbol{S}(t) = (s_{ij}(t))_{2n \times 2n}$，$g(t)$为$2n$维向量值函数，表达式如下：

$$\begin{cases} g_{2i-1} = -\frac{1}{T\eta}\sum_{j=1}^{n} z_j \frac{\partial}{\partial y_i}\left(\frac{\partial E}{\partial y_j}\right), \\ g_{2i} = -\frac{1}{T\eta}\left(\frac{\partial E}{\partial y_i} + \sum_{j=1}^{n} z_j \frac{\partial}{\partial z_i}\left(\frac{\partial E}{\partial y_j}\right)\right), \end{cases} \quad (16)$$

$$\begin{cases} s_{2i-1,2j-1} = \sum_h \left(\frac{\partial a_{jh}}{\partial y_i}\frac{\partial E}{\partial y_j} + a_{jh}\frac{\partial}{\partial y_i}\left(\frac{\partial E}{\partial y_j}\right) \right. \\ \qquad\qquad\left. + \frac{\partial b_{jh}}{\partial y_i}\frac{\partial E}{\partial z_h} + b_{jh}\frac{\partial}{\partial y_i}\left(\frac{\partial E}{\partial z_h}\right)\right), \\ s_{2i-1,2j} = \sum_h \left(\frac{\partial c_{jh}}{\partial y_i}\frac{\partial E}{\partial y_j} + c_{jh}\frac{\partial}{\partial y_i}\left(\frac{\partial E}{\partial y_j}\right) \right. \\ \qquad\qquad\left. + \frac{\partial d_{jh}}{\partial y_i}\frac{\partial E}{\partial z_h} + d_{jh}\frac{\partial}{\partial y_i}\left(\frac{\partial E}{\partial z_h}\right)\right), \\ s_{2i,2j-1} = \sum_h \left(\frac{\partial a_{jh}}{\partial z_i}\frac{\partial E}{\partial y_j} + a_{jh}\frac{\partial}{\partial z_i}\left(\frac{\partial E}{\partial y_j}\right) \right. \\ \qquad\qquad\left. + \frac{\partial b_{jh}}{\partial z_i}\frac{\partial E}{\partial z_h} + b_{jh}\frac{\partial}{\partial z_i}\left(\frac{\partial E}{\partial z_h}\right)\right), \\ s_{2i,2j} = \sum_h \left(\frac{\partial c_{jh}}{\partial z_i}\frac{\partial E}{\partial y_j} + c_{jh}\frac{\partial}{\partial z_i}\left(\frac{\partial E}{\partial y_j}\right) \right. \\ \qquad\qquad\left. + \frac{\partial d_{jh}}{\partial z_i}\frac{\partial E}{\partial z_h} + d_{jh}\frac{\partial}{\partial z_i}\left(\frac{\partial E}{\partial z_h}\right)\right). \end{cases} \quad (17)$$

式(13)是一个变系数的非齐次一阶线性微分方程组的终值问题，它的解具有如下的显示表达式：

$$\lambda(t) = e^{\int_t^{t_f} M(s)ds} \int_t^{t_f} g(s) e^{-\int_s^{t_f} M(v)dv} ds \,. \quad (18)$$

**命题 1**：由极大值原理，最优控制$u^*$及对应最优轨线$x^*$应当满足[19]：

1) $$\begin{cases} \dot{x}_i^*(t) = \frac{\partial H}{\partial \lambda}\big(x_i^*(t), u^*(t), \lambda^*(t)\big), \; i=1,\cdots,n, \\ x^*(t_0) = x_0, \\ \dot{\lambda}^*(t) = -\frac{\partial H}{\partial x}\big(x^*(t), u^*(t), \lambda^*(t)\big) \\ \qquad = \big(S(t)\lambda^*(t) + g(t)\big)\big|_{x^*(t),u^*(t)}, \\ \lambda^*(t_f) = 0, \end{cases}$$

其中$H$为 Halmiton 函数；

2) 使$u^*$连续的所有时间点$t \in (t_0, t_f)$都满足：

$$\frac{\partial H}{\partial u}\big(x^*(t), u^*(t), \lambda^*(t)\big) = 0,$$

$$\frac{\partial^2 H}{\partial u^2}\big(x^*(t), u^*(t), \lambda^*(t)\big) \le 0;$$

3) $H\big(x^*(t_f), u^*(t_f), \lambda^*(t_f)\big) = $ 常数。

## 4．最优泵送性能的数值解法

### 4.1 算法流程

---

**输入**：初始控制输入$u^0(t)$，初始状态$x(t_0)$，算法终止参数$\varepsilon > 0$。

步 0．初始化迭代步$k = 0$。
步 1．解初始值问题

$$\begin{cases} \dot{x}(t) = f\big(x(t), u^k(t)\big), \\ x(t_0) = x_0, \end{cases}$$

并设其解为$x^k(t)$。
步 2．计算$J(u^k)$。
步 3．解终值问题

$$\begin{cases} \dot{\lambda}(t) = S(t)\lambda(t) + g(t), \\ \lambda(t_f) = 0, \end{cases}$$

得$\lambda^k(t)$。
步 4．计算$\nabla J(u^k) = \frac{\partial H}{\partial u}(x^k, u^k, \lambda^k)$。
步 5．计算下降方向系数$\beta^k$。当$k = 0$时，取$\beta^k = 0$；当$k > 0$时，取$\beta^k = \frac{\|\nabla J(u^k)\|^2}{\|\nabla J(u^{k-1})\|^2}$。
步 6．下降方向$p^k = -\nabla J(u^k) + \beta^k p^{k-1}$。
步 7．计算$J(u^{k+1})$，采用可接受算法确定第$k$步的迭代步长因子$\alpha^k$，即$\alpha^k$只须满足$J(u^{k+1}) < J(u^k)$。
步 9．$u^{k+1} = u^k + \alpha^k p^k$。
步 10．验证终止法则$\|\nabla J(u^k)\| < \varepsilon$，满足，则停止计算；否则令$k \Leftarrow k + 1$，转步 1。

**输出**：最优控制$u^* \Leftarrow u^k$，最优泵送值$J(u^*) \Leftarrow J(u^k)$。

---

对于无约束优化问题，常用的数值方法有最速下降法、共轭梯度法、牛顿法、拟牛顿法等。梯度法虽然简单但收敛速度慢，容易在最优附近出现"拉锯"现象；由于牛顿法使用了目标函数的二阶导数信息，具有很快的收敛速度，但其缺点也是十分明显的：求取二阶导数不仅增大了计算工作量及复杂性，也限制了应用范围；拟牛顿法借鉴了牛顿法的思想，同时避免了计算二阶导数，应用范围比牛顿法更广；共轭梯度的梯度方向相互共轭，对于非线性问题近似具有二次终止性质，对比梯度法来说，收敛速度快，对比牛顿法来说，不需要计算二阶导数矩阵。因此本文考虑使用 Fletcher-Reeves 共轭梯度法，给出解决最优控制问题(12)-(13)的数值方法，算法流程如上所示。

### 4.2 仿真结果

我们使用共轭梯度法求解泵送性能最优控制问题(12)-

(13)。具体而言，我们选择初值$u^0(t) \equiv \frac{\pi}{2}$，算法终止参数 $\varepsilon = 10^{-5}$，$x_0 = (0.3a, 0.6a, \cdots, 0.3na)^T$，其余参数设置记录在表 1 中。

表 1 参数设置

| 参数 | 取值 |
| --- | --- |
| $n$ | 20 |
| $\chi$ | 1.704 |
| $a(\mu m)$ | 0.5 |
| $\eta(Ns/m^2)$ | $10^{-3}$ |
| $k(N/m)$ | $1.5 \times 10^{-9}$ |
| $A(Nm)$ | $4.5 \times 10^{-22}$ |
| $B(T)$ | 0.07 |

数值实验结果表明算法迭代 4 步终止。下图 2 记录了磁场控制输入$u(t)$在各迭代步的具体取值。

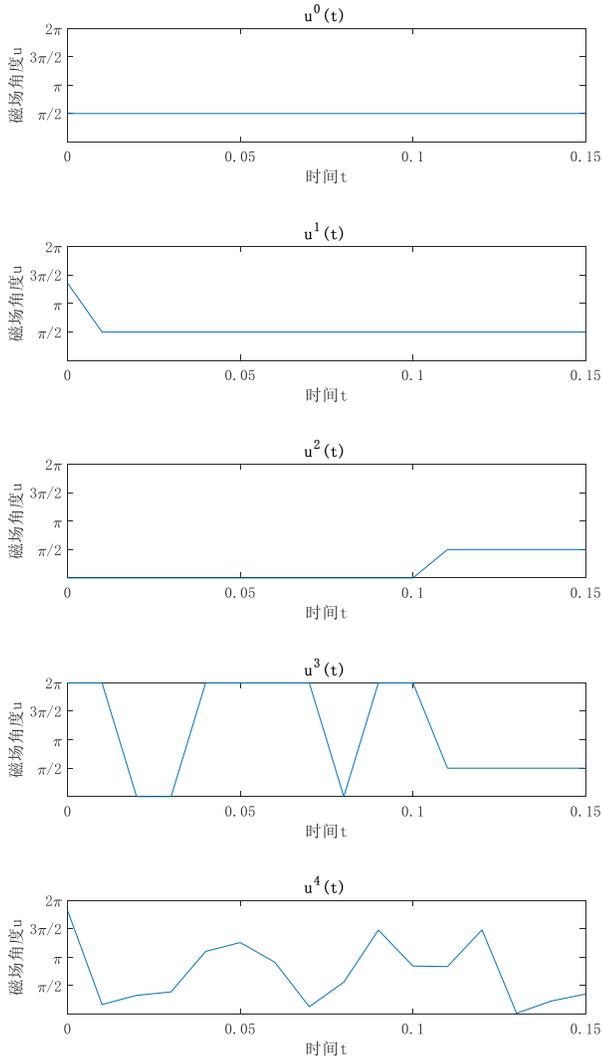

图 2 控制输入$u(t)$在各迭代步的取值

此外，在数值实验中还发现泵送性能随迭代步数增加而快速提高，它们的具体取值如下表 2 所示。进一步，我们还将求解得到的泵送性能$J(u^*)$与其它文献作对比。在同样的参数条件下，求得文献[17]中的磁场控制设置所对应的泵送性能为 $1.14 \times 10^{-13}$Nm，这比我们得到的数值结果 $4.78 \times 10^{-11}$Nm 要小 2 个数量级。

表 2 历次迭代下的泵送性能

| 迭代次数 | 泵送性能(Nm) |
| --- | --- |
| 0 | $1.36 \times 10^{-14}$ |
| 1 | $1.11 \times 10^{-13}$ |
| 2 | $5.86 \times 10^{-13}$ |
| 3 | $1.22 \times 10^{-11}$ |
| 4 | $4.78 \times 10^{-11}$ |

由于共轭梯度法对初值较为敏感，下面我们将选择不同的初值$u^0(t)$进行数值实验。其中$u_1^0(t), u_2^0(t), u_3^0(t), u_4^0(t)$分别为$[0,2\pi]$上不同的随机数函数，$u_5^0(t) \equiv \frac{\pi}{3}$，$u_6^0(t) = -\left(t - \frac{\pi}{4}\right)^2, u_7^0(t) = 3t, u_8^0(t) \equiv 0$，不同初值下算法得到的解$u^*(t)$和$J(u^*)$汇总到下图表。

表 3 不同初值下最优控制问题的数值解结果

| 初值 | 迭代次数 | $\|\nabla J(u^*)\|$ | $J(u^*)$(Nm) |
| --- | --- | --- | --- |
| $u_1^0$ | 4 | 0.0011 | $3.27 \times 10^{-11}$ |
| $u_2^0$ | 8 | 0.00095 | $2.71 \times 10^{-11}$ |
| $u_3^0$ | 8 | 0.0013 | $7.50 \times 10^{-13}$ |
| $u_4^0$ | 9 | 0.066 | $1.53 \times 10^{-10}$ |
| $u_5^0$ | 8 | 0.0059 | $9.07 \times 10^{-12}$ |
| $u_6^0$ | 11 | 0.025 | $1.14 \times 10^{-12}$ |
| $u_7^0$ | 4 | 0.094 | $4.61 \times 10^{-13}$ |
| $u_8^0$ | 2 | 0.026 | $1.92 \times 10^{-13}$ |
| $u_9^0$ | 3 | 0.00018 | $9.98 \times 10^{-12}$ |
| $u_{10}^0$ | 2 | 0.00067 | $1.09 \times 10^{-11}$ |

从数值实验结果可以看出，在我们的模型和最优控制算法下求得的泵送性能$J(u^*)$要优于文献[17]。根据表3，我们将表现不错的初值$u_1^0(t), u_2^0(t), u_4^0(t)$和$u_{10}^0(t)$对应的数值解汇总到下图中3。其中$t_0 = \frac{\pi}{10}, T = \frac{\pi}{10}$。

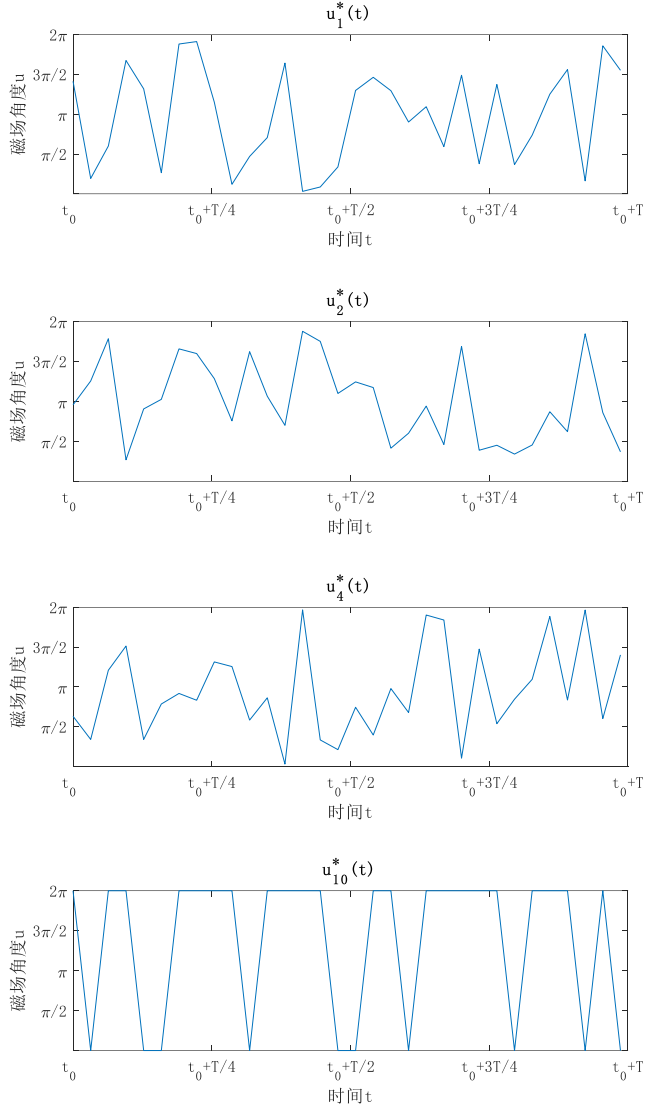

图3 $u_1^*(t), u_2^*(t), u_4^*(t)$和$u_{10}^*(t)$

由于问题(12)-(13)是一个非线性非凸最优控制问题，因此上述算法求得的解$u^*$可能是局部最优的。后面可以考虑使用凸松弛等手段将此非凸问题转化为凸优化问题，进而得到更好的次优解。

## 5．结论与展望

自 Honda 等[20]于 1996 年提出第一个纤毛机器人以来，已经吸引了大量的研究兴趣。然而，如何构建一个实用、高效又精确的纤毛机器人运动物理模型是目前学者们面临的共同难题。近年来多个体系统理论引起了人们的极大兴趣，并被应用到许多领域，包括物理、生物、化学、社会和工程[21-25]，然而在微纳机器人领域应用还很少。

为了探讨如何设计磁场控制方法来诱导纤毛产生最强泵送量，本文将单根人工纤毛建模成一个受磁场驱动的多个体系统，并给出了输送微流体能力的衡量指标—泵送性能及其相应的最优控制问题。此外，利用共轭梯度法对最优控制问题给出了数值解法，并将求解结果与文献[17]进行了对比。仿真结果验证了本文所提出的建模方法以及最优控制算法的有效性。

在实际应用中，为了提高效率一般采用多纤毛系统对微流体进行输送。因此，一个自然的问题是如何对多纤毛系统的泵送性能进行最优控制。然而，多纤毛间的流体动力相互作用非常复杂，并且其协同运动与流体输送的关联尚在研究中，还未有成熟的理论，这给多纤毛系统的建模与优化带来了挑战，需要未来作进一步研究。